\documentclass[12pt]{article}
\usepackage[Lenny]{fncychap}
\usepackage{amsthm}
\usepackage{graphicx} 
\usepackage{tikz}
\usepackage{mathrsfs}
\usepackage{fancyhdr}
\usetikzlibrary{decorations.fractals}
\usepackage[english]{babel}
\usepackage[utf8]{inputenc}
\usepackage[many]{tcolorbox}
\usepackage[a4paper, total={6.5in,10in}]{geometry}
\usepackage{amsfonts,amsmath}
\usepackage{amsmath}
\usepackage{tkz-euclide}
\usepackage{hyperref}
\hypersetup{colorlinks=true,linkcolor=blue,filecolor=magenta,urlcolor=cyan}
\usepackage{float}
\usepackage{lscape}
\usepackage{setspace}

\newtheorem{example}{Example}[subsection]
\newtheorem{coro}{Corollary}[subsection]
\newtheorem{theorem}{Theorem}[section]

\usepackage{mathtools}

\title{\textbf{Series Associated with Harmonic Numbers, Fibonacci Numbers and Central Binomial Coefficients $\binom{2n}{n}$}}
\author{Akerele Olofin Segun$^1$ and Salami Olamide Esther$^2$
\\ Department of Mathematics, University of Ibadan, Oyo State, Nigeria.\\
$^{1}$email: akereleolofin@gmail.com\\
$^{2}$email: salamiolamideesther03@gmail.com}
\date{}
\AtBeginDocument{}

\begin{document}
	\maketitle
\fancyhead{}
\fancyhead[RO]{A.S. Olofin and S.E. Olamide} 
\begin{abstract}
	We find various series that involves the central binomial coefficients $\binom{2n}{n}$, harmonic numbers and Fibonacci Numbers.\\
Contrary to the traditional hypergeometric function $_pF_q$ approach, our method utilizes a straightforward transformation to obtain new evaluations linked to Fibonacci numbers and the golden ratio. Before the end of this paper, we also gave a new series representation for $\zeta(2)$.

\end{abstract}
\section{Introduction}	
In an article published in the year 2016, Hongwei [1], gave a generating function for the sequence, $\binom{2n}{n}H_n$, $\binom{2n}{n}(H_{2n}-H_{n})$, $C_n(H_{2n}-H_n)$ and few others, where $C_n$ is the nth Catalan number, $\displaystyle{C_n=\frac{1}{n+1}\binom{2n}{n}}$. In search of interesting series associated with central binomial coefficients and Harmonic numbers, Hongwei established several interesting sums as follows;
\begin{equation}
	\label{1}
	\sum_{n=1}^{\infty}\frac{1}{4^n(2n+1)}\binom{2n}{n}(H_{2n}-H_n)=\pi\ln 2 - 2G,
\end{equation}
\begin{equation}
	\label{2}
	\sum_{n=1}^{\infty}\frac{1}{4^nn(2n+1)}\binom{2n}{n}(H_{2n-1}-H_n)=2+2\ln 2 + \ln^2 2 + 4G - \pi(1+2\ln 2),
\end{equation}	
\begin{equation}
	\label{3}
	\sum_{n=1}^{\infty}\frac{1}{4^n(2n+3)}C_nH_n=2+4\ln 2 -4G -\pi + \pi \ln 2,
\end{equation}
where $G$ is the Catalan's constant, which is defined by 
$$G:= \sum_{k=0}^{\infty}\frac{(-1)^k}{(2k+1)^2}.$$	
From [1], we have the generating function;
\begin{equation}
	\label{4}
\mathscr{M}(x):=\sum_{n=1}^{\infty}\binom{2n}{n}H_n x^n =\frac{2}{\sqrt{1-4x}}\ln \left(\frac{1+\sqrt{1-4x}}{2\sqrt{1-4x}}\right)
\end{equation}	
which converges on $[-1/4,1/4)$. \\
Now before we continue, let $F_n$ and $L_n$ denote the n-th Fibonacci and Lucas numbers, both satisfying the recurrence relation $\Gamma_n=\Gamma_{n-1}+\Gamma_{n-2}$, $n\geq 2$, with conditions $F_0=0,F_1=1$ and $L_0=2,L_1=1$. Also $L_{-m}=(-1)^mL_m$ and $F_{-m}=(-1)^{m-1}F_m$. Through out this paper, we denote the golden ratio $\alpha=\frac{1+\sqrt{5}}{2}$ and write it's conjugate $\beta=\frac{1-\sqrt{5}}{2}$, so that $\alpha \beta = -1$ and $\alpha + \beta = 1$. We have the Binet formulas for Fibonacci and Lucas numbers to be; 
$$F_m = \frac{\alpha^m-\beta^m}{\alpha-\beta}, \quad L_m = \alpha^m + \beta^m$$
for any integer $m$.\\
We will require the following, which are consequences of the Binet formula and well-known identities which are valid for integers $m$ and $n$.
$$\alpha^{2m}=\alpha^{m}F_m\sqrt{5}-(-1)^{m+1},$$
$$\alpha^{2m}=\alpha^{m}L_m-(-1)^{m},$$
$$\beta^{2m}=\beta^{m}L_{m}-(-1)^m,$$
$$F_n^2+(-1)^{n+m-1}F_m^2=F_{n-m}F_{n+m},$$
$$L_{n+m}+(-1)^mL_{n-m}=L_nL_m.$$
Now setting $\displaystyle{x=\frac{1}{4(\alpha^{2r}+(-1)^{r+1})}}$ in (\ref{4}) for $r\in \mathbb{N}$ we get that;
\begin{equation}
	\label{5}
	\sum_{n=1}^{\infty}\frac{1}{(4\sqrt{5})^n\alpha^{rn}F_r^n}\binom{2n}{n}H_n=\frac{2\sqrt{\alpha^rF_r\sqrt{5}}}{\sqrt{\alpha^rF_r\sqrt{5}-1}}\ln \left(\frac{\sqrt{\alpha^rF_r\sqrt{5}}+\sqrt{\alpha^rF_r\sqrt{5}-1}}{2\sqrt{\alpha^rF_r\sqrt{5}-1}}\right)
\end{equation}	
Evaluation at $r=1,2,3$ in (\ref{5}) gives;
\begin{equation}
	\label{6}
		\sum_{n=1}^{\infty}\frac{1}{(2^2\sqrt{5})^n\alpha^{n}}\binom{2n}{n}H_n = \frac{2\sqrt{\alpha\sqrt{5}}}{\sqrt{\alpha\sqrt{5}-1}}\ln \left(\frac{\sqrt{\alpha\sqrt{5}}+\sqrt{\alpha\sqrt{5}-1}}{2\sqrt{\alpha\sqrt{5}-1}}\right)
\end{equation}	
\begin{equation}
	\label{7}
		\sum_{n=1}^{\infty}\frac{1}{(2^2\sqrt{5})^n\alpha^{2n}}\binom{2n}{n}H_n = \frac{2\alpha\sqrt[4]{5}}{\sqrt{\alpha^2\sqrt{5}-1}}\ln \left(\frac{\alpha\sqrt[4]{5}+\sqrt{\alpha^2\sqrt{5}-1}}{2\sqrt{\alpha^2\sqrt{5}-1}}\right)
\end{equation}
\begin{equation}
	\label{8}
		\sum_{n=1}^{\infty}\frac{1}{(2^3\sqrt{5})^n\alpha^{3n}}\binom{2n}{n}H_n = 	 \frac{2\sqrt{2}\sqrt{\alpha^3\sqrt{5}}}{\sqrt{2\alpha^3\sqrt{5}-1}}\ln \left(\frac{\sqrt{2\alpha^3\sqrt{5}}+\sqrt{2\alpha^3\sqrt{5}-1}}{2\sqrt{2\alpha^3\sqrt{5}-1}}\right)
\end{equation}
Also from (\ref{5}), by replacing $r$ with $2r$ we get,
\begin{equation}
	\label{9}
	\sum_{n=1}^{\infty}\frac{1}{(4\sqrt{5})^n\alpha^{2rn}F_{2r}^n}\binom{2n}{n}H_n=\frac{2\alpha^r\sqrt{F_{2r}\sqrt{5}}}{\sqrt{\alpha^{2r}F_{2r}\sqrt{5}-1}}\ln \left(\frac{\alpha^r\sqrt{F_{2r}\sqrt{5}}+\sqrt{\alpha^{2r}F_{2r}\sqrt{5}-1}}{2\sqrt{\alpha^{2r}F_{2r}\sqrt{5}-1}}\right)
\end{equation}	
In this paper, by playing around the $\mathscr{M}(x)$ in (\ref{4}), we shall produce more interesting results. To ensure accuracy, all formulas appearing in this paper were verified numerically by Mathematica 13.3.
\section{Main Theorems}	
\begin{theorem}
	If r is a natural number, then
	\begin{equation}
		\label{10}
		\sum_{n=1}^{\infty}\frac{1}{4^n\alpha^{rn}L^n_r}\binom{2n}{n}H_n = \frac{2\sqrt{\alpha^r L_r}}{\sqrt{\alpha^rL_r-1}}\ln \left(\frac{\sqrt{\alpha^r L_r}+\sqrt{\alpha^rL_r-1}}{2\sqrt{\alpha^rL_r-1}}\right)
	\end{equation}
\end{theorem}	
\begin{proof}
	Setting $\displaystyle{x=\frac{1}{4(\alpha^{2r}+(-1)^r)}}$ in (\ref{4}) and using the second identity from the list of consequences from the Binet's formula (\emph{introduction}). The result follows immediately.
\end{proof}	
\begin{example}
	Evaluation at $r=1,2,3$ in (\ref{10}), gives;
	\begin{equation}
		\label{11}
		\sum_{n=1}^{\infty}\frac{1}{4^n\alpha^n}\binom{2n}{n}H_n = \frac{2\sqrt{\alpha}}{\sqrt{\alpha-1}}\ln \left(\frac{\sqrt{\alpha}+\sqrt{\alpha -1}}{2\sqrt{\alpha -1}}\right),
	\end{equation}
	\begin{equation}
		\label{12}
		\sum_{n=1}^{\infty}\frac{1}{12^n\alpha^{2n}}\binom{2n}{n}H_n=\frac{2\alpha\sqrt{3}}{\sqrt{3\alpha^2-1}}\ln \left(\frac{\alpha\sqrt{3}+\sqrt{3\alpha^2 -1}}{2\sqrt{3\alpha^2 -1}}\right),	
	\end{equation}
		\begin{equation}
		\label{13}
		\sum_{n=1}^{\infty}\frac{1}{16^n\alpha^{3n}}\binom{2n}{n}H_n=\frac{4\sqrt{\alpha^3}}{\sqrt{4\alpha^3-1}}\ln \left(\frac{2\sqrt{\alpha^3}+\sqrt{4\alpha^3 -1}}{2\sqrt{4\alpha^3 -1}}\right).
	\end{equation}
\end{example}	
\begin{coro}
If r is a natural number, then
\begin{equation}
	\label{14}
	\sum_{n=1}^{\infty}\frac{1}{4^n\alpha^{2rn}L^n_{2r}}\binom{2n}{n}H_n = \frac{2\alpha^r\sqrt{L_{2r}}}{\sqrt{\alpha^{2r}L_{2r}-1}}\ln \left(\frac{\alpha^r\sqrt{ L_{2r}}+\sqrt{\alpha^{2r}L_{2r}-1}}{2\sqrt{\alpha^{2r}L_{2r}-1}}\right).
\end{equation}
\end{coro}	
\begin{proof}
	Replace $r$ with $2r$ in (\ref{10})
\end{proof}	
\newpage
\begin{theorem}
	If $r$ is a non-negative integer, then 
	\begin{equation}
		\label{15}
		\sum_{n=1}^{\infty}\frac{1}{4^n\alpha^{rn}L^n_r}\binom{2n}{n}(H_{2n}-H_n)=-\frac{\sqrt{\alpha^rL_r}}{\sqrt{\alpha^r L_r-1}}\ln \left(\frac{\sqrt{\alpha^rL_r}+\sqrt{\alpha^r L_r-1}}{2\sqrt{\alpha^r L_r}}\right), \quad r>0	
	\end{equation}
\begin{equation}
	\label{16}
	\sum_{n=1}^{\infty}\frac{1}{(4\sqrt{5})^n\alpha^{rn}F^n_r}\binom{2n}{n}(H_{2n}-H_n)=-\frac{\sqrt{\alpha^rF_r\sqrt{5}}}{\sqrt{\alpha^r F_r \sqrt{5}-1}}\ln\left(\frac{\sqrt{\alpha^rF_r\sqrt{5}}+\sqrt{\alpha^rF_r\sqrt{5}-1}}{2\sqrt{\alpha^rF_r\sqrt{5}}}\right), \quad r>0	
\end{equation}
\end{theorem}
\begin{proof}
	Before we establish the following identities above, we shall establish a generating function for the sequence $\binom{2n}{n}(H_{2n}-H_n)$ in the ways of [1]. Now observe, 
	\begin{align*}
	\sum_{n=1}^{\infty}\binom{2n}{n}(H_{2n}-H_{n})x^n = \sum_{n=1}^{\infty}\sum_{k=1}^n\frac{1}{k}\binom{2n-k}{n}x^n=\sum_{k=1}^{\infty}\frac{1}{k}\left\{\sum_{n=k}^{\infty}\binom{2n-k}{n}x^n\right\}\\
\end{align*}

By setting $(n=m+k)$, we have that:
\begin{align*}
\sum_{k=1}^{\infty}\frac{1}{k}\left\{\sum_{n=k}^{\infty}\binom{2n-k}{n}x^n\right\}=	\sum_{k=1}^{\infty}\sum_{m=0}^{\infty}\binom{2m+k}{m+k}x^{m+k}\\
=\sum_{k=1}^{\infty}\frac{x^k}{k}\left\{\sum_{m=0}^{\infty}\binom{2m+k}{m}x^m\right\}
\end{align*}
Since, (from [2, p.203])
$$ \sum_{m=0}^{\infty}\binom{2m+k}{m}\vartheta^m= \frac{1}{\sqrt{1-4\vartheta}}\left(\frac{1-\sqrt{1-4\vartheta}}{2\vartheta}\right)^k, \quad |\vartheta|<1/4  $$
Thus, we have that 
\begin{equation}
\label{17}
	\sum_{n=1}^{\infty}\binom{2n}{n}(H_{2n}-H_{n})x^n =\frac{1}{\sqrt{1-4x}}\sum_{k=1}^{\infty}\frac{1}{k}\left(\frac{1-\sqrt{1-4x}}{2}\right)^k =-\frac{1}{\sqrt{1-4x}}\ln\left(\frac{1-\sqrt{1-4x}}{2}\right)
\end{equation}
Now set $\displaystyle{x=\frac{1}{4(\alpha^{2r}+(-1)^{r+1})}}$ and $\displaystyle{x=\frac{1}{4(\alpha^{2r}+(-1)^{r})}}$ in (\ref{17}), then (\ref{16}) and (\ref{15}) follows directly.
\end{proof}
\begin{example}
	Evaluation at $r=1, 2, 3$ in (\ref{15}) and (\ref{16}), respectively, gives
\begin{equation}
	\label{18}
\sum_{n=1}^{\infty}\frac{1}{4^n\alpha^n}\binom{2n}{n}(H_{2n}-H_{n})=-\frac{\sqrt{\alpha}}{\sqrt{\alpha-1}}\ln\left(\frac{\sqrt{\alpha}+\sqrt{\alpha -1}}{2\sqrt{\alpha}}\right)
\end{equation}	
\begin{equation}
	\label{19}
\sum_{n=1}^{\infty}\frac{1}{12^n\alpha^{2n}}\binom{2n}{n}(H_{2n}-H_{n})=-\frac{\alpha\sqrt{3}}{\sqrt{3\alpha^2-1}}\ln\left(\frac{\alpha\sqrt{3}+\sqrt{3\alpha^2 -1}}{2\alpha\sqrt{3}}\right)
\end{equation}	
\newpage
\begin{equation}
	\label{20}
\sum_{n=1}^{\infty}\frac{1}{16^n\alpha^{3n}}\binom{2n}{n}(H_{2n}-H_{n})=-\frac{2\sqrt{\alpha^3}}{\sqrt{4\alpha^3-1}}\ln\left(\frac{2\sqrt{\alpha^3}+\sqrt{4\alpha^3 -1}}{4\sqrt{\alpha^3}}\right)
\end{equation}	
\begin{equation}
	\label{21}
\sum_{n=1}^{\infty}\frac{1}{(4\sqrt{5})^n\alpha^{n}}\binom{2n}{n}(H_{2n}-H_{n})=-\frac{\sqrt{\alpha\sqrt{5}}}{\sqrt{\alpha\sqrt{5}-1}}\ln\left(\frac{\sqrt{\alpha\sqrt{5}}+\sqrt{\alpha\sqrt{5}-1}}{2\sqrt{\alpha\sqrt{5}}}\right)
\end{equation}	
\begin{equation}
	\label{22}
\sum_{n=1}^{\infty}\frac{1}{(4\sqrt{5})^n\alpha^{2n}}\binom{2n}{n}(H_{2n}-H_{n})=-\frac{\alpha\sqrt[4]{5}}{\sqrt{\alpha^2\sqrt{5}-1}}\ln \left(\frac{\alpha \sqrt[4]{5}+\sqrt{\alpha^2\sqrt{5}-1}}{2\alpha\sqrt[4]{5}}\right)
\end{equation}
\begin{equation}
	\label{23}
\sum_{n=1}^{\infty}\frac{1}{(8\sqrt{5})^n\alpha^{3n}}\binom{2n}{n}(H_{2n}-H_{n})=-\frac{\sqrt{2}\sqrt{\alpha^3\sqrt{5}}}{\sqrt{\alpha^32\sqrt{5}-1}}\ln\left(\frac{\sqrt{2}\sqrt{\alpha^3\sqrt{5}}+\sqrt{\alpha^32\sqrt{5}-1}}{2\sqrt{2}\sqrt{\alpha^3\sqrt{5}}}\right)
\end{equation}	
\end{example}
\begin{theorem}
$$\sum_{n=1}^{\infty}\frac{C_n(H_{2n}-\frac{1}{2}H_n)}{4^n(2n+1)}\left(\frac{\pi}{2}-\frac{(2n)!!}{(2n+1)!!}\right)=2\ln(2)+\frac{7}{8}\zeta(3)+\frac{\pi}{12}(-12+\pi(-1+\ln(8)))$$
where $C_n$ is the n-th Catalan's number and , 
$$n!! := \prod_{k=0}^{\lceil\frac{n}{2}\rceil-1}(n-2k):=\begin{cases}
	\vspace{3mm}
\displaystyle{\prod_{k=1}^{\frac{n}{2}}(2k)}, \text{if $n$ is even},\\
\displaystyle{\prod_{k=1}^{\frac{n+1}{2}}(2k-1)}, \text{if $n$ is odd}.
\end{cases}$$
\end{theorem}	
\begin{proof}
	From [1] we have that;
$$\sum_{n=1}^{\infty}\binom{2n}{n}(H_{2n}-\frac{1}{2}H_n)x^n = -\frac{1}{\sqrt{1-4x}}\ln \sqrt{1-4x} $$
Immediately we see that;
$$\sum_{n=1}^{\infty}C_n(H_{2n}-\frac{1}{2}H_n)x^n = \frac{1}{2x}(1-\sqrt{1-4x}+\sqrt{1-4x}\ln \sqrt{1-4x}), \quad |x|<1/4$$
For the interesting part, set $\displaystyle{x=\frac{1}{4}\sin^2t}$ for $t\in (-\pi/2,\pi/2)$. Then we have that;
\begin{equation}
	\label{24}
\sum_{n=1}^{\infty}\frac{C_n(H_{2n}-\frac{1}{2}H_n)}{4^n}\sin^{2n}t=\frac{2}{\sin^2t}(1-\cos t + \cos t\ln \cos t).
\end{equation}
Now multiply (\ref{24}) by $t\cos t$ and integrating both sides from 0 to $\pi/2$. We have;
\begin{align*}
	\sum_{n=1}^{\infty}\frac{C_n(H_{2n}-\frac{1}{2}H_n)}{4^n}\int_{0}^{\frac{\pi}{2}}t\cos t\sin^{2n}t\, dt = 2\int_{0}^{\frac{\pi}{2}}\frac{t\cos t}{\sin^2t}(1-\cos t+\cos t \ln \cos t)\,dt
\end{align*}
Thus, we evaluated the integral on the right hand side to $\displaystyle{\ln 2 + \frac{\pi}{24}(-12+\pi(-1+\ln 8))+\frac{7}{16}\zeta(3)}$ and using integration by part for the integral on the left we have that; 
\begin{equation*}
\int_{0}^{\frac{\pi}{2}}\frac{t\cos t}{\sin^2t}(1-\cos t+\cos t \ln \cos t)\,dt=\ln 2 + \frac{\pi}{24}(-12+\pi(-1+\ln 8))+\frac{7}{16}\zeta(3),	
\end{equation*}
\begin{equation*}
\int_{0}^{\frac{\pi}{2}}t\cos t\sin^{2n}t\, dt =\frac{1}{2n+1}\left(\frac{\pi}{2}-\frac{(2n)!!}{(2n+1)!!}\right).
\end{equation*}
So the result follows immediately.
\end{proof}	

\begin{theorem}(\textbf{\textit{Two Ramanujan-Like Series}})
	$$\sum_{n=1}^{\infty}\frac{C_n(H_{2n}-H_n)}{4^{2n}}\binom{2n+2}{n+1}=\frac{16}{\pi}\psi,$$
where $\psi = 2G+\pi-2-\ln 2-\pi \ln 2$.
$$\sum_{n=1}^{\infty} \frac{C_nH_{2n}}{16^n}\binom{2n}{n}=\frac{2}{\pi}\psi^*,$$
where $\psi^*= 2+\pi -2\ln 8.$
\end{theorem}	
\begin{proof}
Integrating both sides of (\ref{17}) with respect to $x$, we obtain the generating function for the sequence $C_n(H_{2n}-H_n)$ as follows;
\begin{equation}
\label{25}
\sum_{n=1}^{\infty}C_n(H_{2n}-H_n)x^n = \frac{1}{2x}\left[(1-\sqrt{1-4x})+(1+\sqrt{1-4x})\ln \left(\frac{1+\sqrt{1-4x}}{2}\right)\right]
\end{equation}
Likewise, from [1] we have that;
\begin{equation}
\label{26}
\sum_{n=1}^{\infty}\binom{2n}{n}H_{2n}x^{n}=\frac{1}{\sqrt{1-4x}}\left[\ln \left(\frac{1+\sqrt{1-4x}}{2}\right)-2\ln \sqrt{1-4x}\right]
\end{equation}
For this reason, integrating both sides of (\ref{26}) we have that,
\begin{equation}
\label{27}
\sum_{n=1}^{\infty}C_nH_{2n}x^n=\frac{1}{2x}[(1-\sqrt{1-4x})-(1+\sqrt{1-4x})\ln (1+\sqrt{1-4x})+\ln 2 + \sqrt{1-4x}\ln (2-8x)]
\end{equation}
Now set $\displaystyle{x=\frac{sin^2t}{4}}$ in (\ref{25}) and (\ref{27}) and integrating both sides with respect to $x$ from 0 to $\pi/2$, while using the following results;
$$\int_{0}^{\frac{\pi}{2}}\left\{(1-\cos t)+(1+\cos t)\ln \left(\frac{1+\cos t}{2}\right)\right\}\, dt = 2G+\pi -2-\ln 2-\pi \ln 2,$$
$$\int_0^{\frac{\pi}{2}}\frac{2}{\sin^2 t}\left[(1-\cos t)-(1+\cos t)\ln (1+\cos t)+\ln 2+\cos t\ln(2\cos^2t)\right]\, dt = 2+\pi -2\ln 8,$$
and, 
$$\int_{0}^{\frac{\pi}{2}}\sin^{2n}t\,dt = \frac{\pi}{2}\frac{(2n)!}{4^n(n!)^2}.$$
Then the series follows directly.
\end{proof}
\begin{theorem}[\textit{A new series representation for $\zeta(2)$}]
$$\zeta(2)=\sum_{n=1}^{\infty}\frac{1024n}{3(2n-1)^2(2n+1)(2n+3)^2}\frac{\binom{2n}{n}}{\binom{2n+2}{n+1}}$$
where $\displaystyle{\zeta(s)=\sum_{n=1}^{\infty}\frac{1}{n^s}}$, provided $\Re (s)>0$.
\end{theorem}	
\begin{proof}
From [5], we have that;
\begin{equation}
\label{28}
\sum_{n=1}^{\infty}\frac{nx^{2n}}{4^n(2n-1)^2(2n+1)}\binom{2n}{n}=\frac{1}{8}\left(\sqrt{1-x^2}+2x\sin^{-1}x-\frac{\sin^{-1}x}{x}\right).
\end{equation}
Next we multiply (\ref{28}) by $x^2$ and integrate both sides with respect to $x$ to get; 
\begin{equation}
\label{29}
\sum_{n=1}^{\infty}\frac{nx^{2n+3}}{4^n(2n-1)^2(2n+1)(2n+3)}\binom{2n}{n}=\frac{(8x^4-8x^2+3)\sin^{-1}x+\sqrt{1-x^2}(6x^3-3x)}{128}
\end{equation}
which converges for $x\in [-1,1]$. Now, we set $x=\sin t$ in (\ref{29}) and integrate over the interval 0 to $\pi/2$. The result follows.
\end{proof}	
\begin{theorem}[A Ramanujan-Like Series involving the ratio $G/\pi$]
$$	\sum_{n=1}^{\infty}\frac{n^2}{16^n(2n-1)^2(2n+1)}\binom{2n}{n}^2 = \frac{G}{4\pi}+\frac{1}{8\pi}$$
\end{theorem}	
\begin{proof} We begin by differentiating (\ref{28}) with respect to $x$ to get,
\begin{equation}
\label{30}
\sum_{n=1}^{\infty}\frac{2n^2x^{2n-1}}{4^n(2n-1)^2(2n+1)}\binom{2n}{n}=\frac{1}{8x^2}\left((2x^2+1)\sin^{-1}x-x\sqrt{1-x^2}\right).
\end{equation}	
Now multiply (\ref{30}) by $x$ and set $x=\sin t$, then we integrate both sides over the interval 0 to $\pi/2$. Finally, using the result; 
$$\int_0^{\frac{\pi}{2}}\frac{t(2\sin^2t+1)-\sin t\cos t}{\sin t}\, dt=\int_0^{\frac{\pi}{2}} (2t\sin t -\cos t)\, dt + \int_0^{\frac{\pi}{2}} \frac{t}{\sin t}\, dt = 1+2G$$
The Series follows immediately.
\footnote{Now we look at some other interesting Series}
\end{proof}
\newpage
\section{Some Deluxe Series}	
Observe that the series in (\ref{25}) converges on $[-1/4,1/4)$. Setting $x=-1/8,1/16$ and $-1/16$ respectively, we obtain the following series;
\begin{equation}
\sum_{n=1}^{\infty}\frac{(-1)^nC_n}{8^n}(H_{2n}-H_n)=-\frac{4}{\sqrt{2}}\left[(\sqrt{2}-\sqrt{3})+(\sqrt{2}+\sqrt{3})\ln\left(\frac{\sqrt{2}+\sqrt{3}}{2\sqrt{2}}\right)\right],
\end{equation}	
\begin{equation}
\sum_{n=1}^{\infty}\frac{C_n}{16^n}(H_{2n}-H_n)=4\left[(2-\sqrt{3})+(2+\sqrt{3})\ln\left(\frac{2+\sqrt{3}}{4}\right)\right],
\end{equation}	
\begin{equation}
\sum_{n=1}^{\infty}\frac{(-1)^nC_n}{16^n}(H_{2n}-H_n)=-4\left[(2-\sqrt{5})+(2+\sqrt{5})\ln\left(\frac{2+\sqrt{5}}{4}\right)\right].
\end{equation}	
In the same light, recall from (\ref{29}) converges on $[-1,1]$. Setting $x=1$ and -1 respectively, we get two interesting series;
\begin{equation}
\sum_{n=1}^{\infty}\frac{n}{4^n(2n-1)^2(2n+1)(2n+3)}\binom{2n}{n}=\frac{3\pi}{256},
\end{equation}
\begin{equation}
\sum_{n=1}^{\infty}\frac{(-1)^{2n+3}n}{4^n(2n-1)^2(2n+1)(2n+3)}\binom{2n}{n}=\frac{-3\pi}{256}.
\end{equation}	
From (\ref{30}) we have;
\begin{equation}
	\sum_{n=1}^{\infty}\frac{n^2}{4^n(2n-1)^2(2n+1)}\binom{2n}{n}=\frac{3\pi}{32}.
\end{equation}	
Also, from (\ref{17}) we get that; 	
\begin{equation}
\sum_{n=1}^{\infty}\frac{H_{2n}-H_n}{8^n}\binom{2n}{n}=-\sqrt{2}\ln \left(\frac{\sqrt{2}-1}{2\sqrt{2}}\right),
\end{equation}	
\begin{equation}
\sum_{n=1}^{\infty}\frac{H_{2n}-H_n}{16^n}\binom{2n}{n}=-\frac{2}{\sqrt{3}}\ln \left(\frac{2-\sqrt{3}}{2\sqrt{2}}\right),
\end{equation}	
\begin{equation}
\sum_{n=1}^{\infty}\frac{C_n}{8^n}\left(H_{2n}-\frac{1}{2}H_n\right)=4\left(1-\frac{1}{\sqrt{2}}-\frac{\ln 2}{2\sqrt{2}}\right),
\end{equation}	
\begin{equation}
\sum_{n=1}^{\infty}\frac{C_n}{16^n}\left(H_{2n}-\frac{1}{2}H_n\right)=8\left(1-\frac{\sqrt{3}}{2}+\frac{\sqrt{3}}{2}\ln \frac{\sqrt{3}}{2}\right).
\end{equation}	
\section{Conclusion}
\begin{itemize}
	\item [1.] To assure accuracy of the results, we verified all the series via \textit{Mathematica 13.3}.
	\item [2.] In this paper we presented new closed forms for some types of series involving the Central Binomial coefficients $\binom{2n}{n}$. To prove our results, we used some generating functions, combined with basic differentiation and integration. Using similar techniques, we  established series evaluations involving Harmonic numbers with Fibonacci and Lucas sequences. In addition, readers can exploit (\ref{28}),(\ref{29}) and (\ref{30}) to generate more exotic series.
\end{itemize}	
\section{Acknowledgment}
We would like to appreciate Prof. U. N. Bassey, Dr. D. A. Dikko and Mr. G. S. Lawal from the Department of Mathematics at the University of Ibadan for their motivations. This helped improve and enrich this article.

\end{document}